\documentclass[leqno,12pt]{article}
\usepackage{inputenc}
\usepackage[english]{babel}
\usepackage{amsmath}
\usepackage{amssymb}
\usepackage{amsthm}
\date{}
\newtheorem{thm}{Theorem}[section]
\newtheorem{lem}[thm]{Lemma}
\newtheorem{prop}[thm]{Proposition}
\newtheorem{cor}[thm]{Corollary}

\theoremstyle{definition}
\newtheorem{ex}[thm]{\it Example}
\newtheorem{rem}[thm]{\it Remark}
\newtheorem{qst}[thm]{\it Question}

\numberwithin{equation}{section}

\title{On the associated graded ring of a semigroup ring}

\begin{document}
\newcommand{\Hi}{\mathrm {Hilb}}
\newcommand{\card}{\mathrm {card}}
\newcommand{\ap}{\mathrm {Ap}}
\newcommand{\ord}{\mathrm {ord}}
\newcommand{\map}{\mathrm {maxAp_M}}
\author{M. D'Anna\thanks{{\em email} mdanna@dmi.unict.it} \and V. Micale \and A. Sammartano}

\maketitle

\begin{abstract}
\noindent Let $(R,\mathbf m)$ be a numerical semigroup ring. In
this paper we study the properties of its associated graded ring
$G(\mathbf m)$. In particular, we describe the $H^0_{\mathcal M}$
for $G(\mathbf m)$ (where $\mathcal M$ is the homogeneous maximal ideal of $G(\mathbf m)$)
and we characterize when $G(\mathbf m)$ is
Buchsbaum. Furthermore, we find the length of $H^0_{\mathcal M}$
as a $G(\mathbf m)$-module, when $G(\mathbf m)$ is Buchsbaum. In
the $3$-generated numerical semigroup case, we describe the
$H^0_{\mathcal M}$ in term of the Apery set of the numerical
semigroup associated to $R$. Finally, we improve two
characterizations of the Cohen-Macaulayness and Gorensteinness of
$G(\mathbf m)$ given in \cite{BF} and \cite{B}, respectively.
\end{abstract}

\section{Introduction}

Let $(R,\mathbf m)$ be a Noetherian, one-dimensional, local ring
with $|R/\mathbf m|=\infty$ and let $G(\mathbf m)=\oplus_{i\geq
0}\mathbf m^i/ \mathbf m^{i+1}$ be the associated graded ring of
$R$ with respect to $\mathbf m$. The study of the properties of
$G(\mathbf m)$ is a classical subject in local algebra.

The concept of a Buchsbaum ring is the most important of all
notions generalizing Cohen-Macaulay rings. While the property for
$G(\mathbf m)$ to be Cohen-Macaulay has been studied extensively
(see, e.g., \cite {BF}, \cite {RV}, or, for the particular case of
semigroup rings, \cite {Ga},\cite {MPT}), not much it is known
about the Buchsbaumness of $G(\mathbf m)$, at least in the general
case (see \cite{Go1} and \cite{Go2}).

As for the semigroup ring case, Sapko, in \cite {S}, gives some
necessary and sufficient conditions for $G(\mathbf m)$ to be
Buchsbaum, when $R$ is associated to a $3$-generated numerical
semigroup; still in the $3$-generated case Shen, in \cite{SH},
studies the Buchsbaumness of $G(\mathbf m)$ and gives positive
answers to the conjectures proposed in \cite{S}. If $S$ is a
general numerical semigroup, it is possible to find some results
on the Buchsbaumness of $G(\mathbf m)$ in \cite{DMM} (where it is
mainly studied the more general case of one dimensional rings) and
in \cite{CZ}.

In this paper we mainly study the Buchsbaumness of $G(\mathbf m)$,
when $(R, \mathbf m)$ is the semigroup ring associated to a
numerical semigroup, but, applying our techniques, we also get
some new results on its Cohen-Macaulayness and on its
Gorensteinness.

In Section 2 we give some preliminaries about numerical semigroups
and semigroup rings associated to a numerical semigroup and we
recall some results on the Buchsbaumness of one dimensional graded
rings proved in \cite{DMM}.

In Section 3 we give a description of $H^0_{\mathcal
M}:=(0:_{G(\mathbf m)}\mathcal M^r)$ (where $r$ is the reduction
number of $\mathbf m$ and where $\mathcal M$ is the homogeneous
maximal ideal of $G(\mathbf m)$ (cf. Corollary \ref{3}) and we use
it in order to characterize when $G(\mathbf m)$ is Buchsbaum (cf.
Proposition \ref{4a}). Successively, we find the length of
$H^0_{\mathcal M}$ as a $G(\mathbf m)$-module when $G(\mathbf m)$
is Buchsbaum (cf. Proposition \ref{10}). Finally, we relate the
Buchsbaumness with a property of the Apery set of the associated
numerical semigroup (cf. Proposition \ref{14}), using a partial
ordering in $S$ introduced in \cite{B}.

In Section 4, we restrict our attention to the semigroup ring
associated to a $3$-generated numerical semigroup $S$; we use the
results of Section 3 in order to prove that, if $G(\mathbf m)$ is
Buchsbaum, then we can determine the $H^0_{\mathcal M}$ in term of
the Apery set of $S$ (cf. Theorem \ref{42}
and Corollary \ref{44}). In particular, we completely solve
\cite[Conjecture 33]{S} and \cite[Conjecture 24]{S}. Finally, we
give a new proof of a result of Shen, which shows that $G(\mathbf
m)$ is Buchsbaum if and only if it is Cohen-Macaulay, for the
$3$-generated symmetric semigroup case (cf. Corollary \ref{45}).

Finally, in Section 5, using the techniques introduced in Section
3, we strengthen, for the semigroup ring case, a characterization
of the Cohen-Macaulayness of $G(\mathbf m)$ given in \cite[Theorem
2.6]{BF} (cf. Proposition \ref{21}). Moreover, we prove that,
assuming the hypotheses of $M$-purity and symmetry for $S$,
$G(\mathbf m)$ is Buchsbaum if and only if it is Cohen-Macaulay
(cf. Proposition \ref{46}) and we use this result to give a
characterization for $G(\mathbf m)$ to be Gorenstein, improving an
analogous result by \cite{B} (cf. Corollary \ref {47}). 
Finally, Question 5.8 about a possible improvement of Proposition \ref{46}
received an affirmative answer by Shen,
that we publish with his permission (cf. Proposition \ref{answer}).

The computations made for this paper are performed by using the GAP
system \cite{GAP} and, in particular, the NumericalSgps package
\cite{NS}.

\section{Preliminaries}

We start this section recalling some well known facts on numerical
semigroups and semigroup rings. For more details see, e.g.,
\cite{BDF}.

A subsemigroup $S$ of the monoid of natural numbers $(\mathbb
N,+)$, such that $0\in S$, is called a {\em numerical semigroup}.
Each numerical semigroup $S$ has a natural partial ordering $\leq
_S$ where, for every $s$ and $t$ in $S$, $s\leq _S t$ if there is
an $u\in S$ such that $t=s+u$. The set $\{g_i\}$ of the minimal
elements in $S\setminus \{0\}$ in this ordering is called {\em the
minimal set of generators} for $S$. In fact all elements of $S$
are linear combinations of minimal elements, with non-negative
integers coefficients. Note that the minimal set $\{g_i\}$ of
generators is finite since for any $s\in S$, $s\neq 0$, we have
that $g_i$ is not congruent to $g_j$ modulo $s$.

A numerical semigroup generated by $g_1<g_2<\dots<g_n$ is denoted
by $\langle g_1, g_2,\dots,g_n\rangle$. Since the semigroup $S=\langle g_1,g_2,\dots,g_n\rangle$ is isomorphic to $\langle
dg_1,dg_2,\dots,dg_n\rangle$ for any $d\in \mathbb N\setminus
\{0\}$, we assume, in the sequel, that $\gcd
(g_1,g_2,\dots,g_n)=1$. It is well known that this condition is
equivalent to $| \mathbb N\setminus S |<\infty$. Hence there is a
well defined the integer $g=g(S)=\max \{x\in \mathbb Z\mid x\notin
S\}$, called {\em Frobenius number} of $S$.

Since the Frobenius number $g$ does not belong to $S$, if $x \in
S$, it is obvious that $g-x\notin S$. A numerical semigroup is
called {\em symmetric} if the converse holds: let $x$ be an
integer, then $g-x \notin S$ implies that $x \in S$.

A {\em relative ideal} of a semigroup $S$ is a nonempty subset $H$
of $\mathbb Z$ such that $H+S\subseteq H$ and $H+s\subseteq S$ for
some $s\in S$. A relative ideal of $S$ which is contained in $S$
is simply called an {\em ideal} of $S$. The ideal $M=\{s\in S \mid
s\neq 0\}$ is called the {\em maximal ideal of $S$}. It is
straightforward to see that, if $H$ and $L$ are relative ideals of
$S$, then $H+L$, $kH(= H+\dots+H$, $k$ summands, for $k\geq 1$)
and $H-_{\mathbb Z}L:=\left\{z\in \mathbb Z\ :\ z+L\subseteq H
\right\}$ are also relative ideals of $S$.

\medskip

The rings $R=k[[t^S]]=k[[t^{g_{1}},\dots,t^{g_{n}}]]$ and
$R=k[t^S]_{\mathbf m}$ are called the {\em numerical semigroup
rings} associated to $S$, where $\mathbf
m=(t^{g_{1}},\dots,t^{g_{n}})$. $R$ is a one-dimensional local
domain with maximal ideal $\mathbf m$ and quotient field
$Q(R)=k((t))$ and $Q(R)=k(t)$, respectively. In both cases the
associated graded ring $G(\mathbf m)$, which is the object of our
investigation, is the same. From now on, we will assume that
$R=k[[t^S]]$, but the other case is perfectly analogous.

We will denote by $v: k((t))\longrightarrow\mathbb Z\cup \infty$
the natural valuation (with associated (discrete) valuation ring
$k[[t]]$), that is
$$
v(\sum_{h=i}^{\infty}r_{h}t^h)=i,\ i\in\mathbb Z,\ r_{i}\ne 0
$$
(in the case $Q(R)=k(t)$, we have the valuation associated to the
DVR $K[t]_{(t)}$). It is straightforward that $v(R)=\{v(r) \mid
r\in R\setminus \{0\}\}=S$.

The relation between $R$ and $S=v(R)$ is very tight and we can
translate many properties of $R$ to the corresponding properties
of $S$. In particular, if $I$ and $J$ are fractional ideals of
$R$, then $v(I)$ and $v(J)$ are relative ideal of $S=v(R)$;
moreover, if $I$ and $J$ are monomial ideals, it is not difficult
to check that $v(I\cap J)=v(I)\cap v(J)$, $v(IJ)=v(I)+v(J)$ and
$v(I:_{Q(R)}J)=v(I)-_{\mathbb Z}v(J)$. Furthermore, if $J\subseteq
I$, then $\lambda_R(I/J)=|v(I)\setminus v(J)|$, where
$\lambda_R(\cdot)$ is the length as $R$-module.

\medskip

Following the notation in \cite{BF}, we denote by
$\ap_{g_1}(S)=\{\omega_0,\dots,\omega_{g_1-1}\}$ the {\em Apery
set} of $S$ with respect of $g_1$, that is, the set of the
smallest elements in $S$ in each congruence class modulo $g_1$.
More precisely, $\omega_0=0$ and $\omega_i=\min\{s\in S\mid
s\equiv i \pmod{g_1}\}$. It is clear that the largest element in
the Apery set is always $g+g_1$. Moreover, if $S$ is symmetric,
then, for every index $j$, there exists an index $i$ such that
$\omega_j+\omega_i=g+g_1$; hence, in the symmetric case, $g+g_1$
is the maximum of the Apery set with respect to $\leq_S$.
Furthermore, it is easy to see that, if $\omega_h+\omega_t\equiv g+g_1$, then
$\omega_h+\omega_t= g+g_1$.

By \cite[Formula I.2.4]{BDF} we have that the blow up of $S$ is
the numerical semigroup $S'=\bigcup_i(iM-_{\mathbb Z}iM)=\langle
g_1,g_2-g_1,\dots,g_n-g_1\rangle$. Note that the set of the
generators $\{g_1,g_2-g_1,\dots,g_n-g_1\}$ is not necessarily the
minimal ones for $S'$; moreover, $g_1$ might not be the smallest
non zero element in $S'$.

In \cite{BF} are defined two families of invariants of $S$, that
give information on the Cohen-Macaulayness of $G(\mathbf m)$. Let
$\ap_{g_1}(S')=\{\omega'_0,\dots,\omega'_{g_1-1}\}$. For each
$i=0,1,\dots,g_1-1$, let $a_i$ be the only integer such that
$\omega'_i+a_ig_1=\omega_i$ and let $b_i=\max \{l\mid\omega_i\in
lM\}$. Clearly $b_0=a_0=0$. Furthermore, $1\le b_i\le a_i$
\cite[Lemma 2.4]{BF}. The following result is proved in a more general setting,
but we give the statement we will need in the sequel, that is for numerical semigroup rings;
notice that under these hypotheses it could be deduced by Remark \ref{1}.

\begin{thm}\label{0}(\cite[Theorem 2.6]{BF})
If $R$ is a semigroup ring, then $G(\mathbf
m)$ is Cohen-Macaulay (briefly, C-M) if and only if \ $a_i=b_i$, for
every $i=0,\dots,g_1-1$.
\end{thm}

A one-dimensional graded ring $T$ with homogeneous maximal ideal
$\mathcal M$ is called {\em Buchsbaum} if $\mathcal M\cdot
H^0_{\mathcal M}=0$ (cf. \cite[Corollary 1.1]{SV}). Since
$H^0_{\mathcal M}=(\cup_{k\ge 1}(0:_T\mathcal M^k))$, the previous
definition is equivalent to
$$
\mathcal M\cdot (\cup_{k\ge 1}(0:_T\mathcal
M^k))=0.
$$

Let $R$ be a Noetherian, one-dimensional, local ring with maximal
ideal $\mathbf m$ such that $|R/\mathbf m|=\infty$ and $\mathbf m$
contains a non-zerodivisor and let $r$ be the reduction number of
$\mathbf m$, that is the minimal natural number such that $\mathbf
m^{r+1}=x\mathbf m^{r}$, with $x$ a superficial element of $R$
(recall that such number $r$ exists by \cite[Theorem 1, Section
2]{NR}). Then we have the following characterization.

\begin{prop}\cite[Corollary 2.3]{DMM}\label{a}
$G(\mathbf m)$ is Buchsbaum if and only if
$(0:_{G(\mathbf m)}\mathcal M)=(0:_{G(\mathbf m)}\mathcal M^{r})$.
\end{prop}

It is also possible to give the graded description of
$(0:_{G(\mathbf m)}\mathcal M^{r})$ as follows (cf. \cite[Formula
(2.3)]{DMM}):

\begin{equation}\label{form1}
\text{$(0:_{G(\mathbf m)}\mathcal M^{r})=\bigoplus_{h=
1}^{r-2}\frac{(\mathbf m^{h+r+1}:_R\mathbf m^{r})\cap\mathbf m^{h}}{\mathbf m^{h+1}}$}.
\end{equation}

Furthermore, if we denote by $R'$ the blow-up of $R$, that is, in
our setting, $R'=\bigcup_i(\mathbf m^i:_{Q(R)}\mathbf
m^i)=R[\frac{\mathbf m}{x}]$ (see e.g. \cite{L}), we have that
$v(R')=S'$ and in \cite[Proposition 2.5]{DMM} it is proved that:

 \begin{equation}\label{form2}
\text{$(0:_{G(\mathbf m)}\mathcal M^{r})=\bigoplus_{h=1}^{r-2}\frac{x^{h+1}R'\cap \mathbf m^{h}}{\mathbf
m^{h+1}}\ .$}
\end{equation}

\begin{rem}\label{1a}
Since the valuation of any superficial element is $v(x)=g_1$, we
can translate as follows the previous formula at the numerical
semigroup level: \noindent let $G(\mathbf m)$ be not C-M and let
$s\in hM\setminus (h+1)M$, then
$$
\overline{t^{s}}\in (0:_{G(\mathbf m)}\mathcal M^r)\
\Longleftrightarrow \ s-(h+1)g_1\in S'.
$$
\end{rem}

\begin{rem}\label{1}
Notice that $a_i>b_i$ if and only if $\overline{t^{\omega_i}}\in (0:_{G(\mathbf
m)}\mathcal M^r)$. Indeed, if $a_i>b_i$, by definition we
have $\omega_i=\omega_i'+a_ig_1= \omega'_i
+(a_i-b_i-1)g_1+(b_i+1)g_1\in b_iM\setminus (b_i+1)M$. Setting
$\alpha=\omega'_i +(a_i-b_i-1)g_1 \in S'$, we get
$\alpha+(b_i+1)g_1=\omega_i$. Hence $\overline{t^{\omega_i}}\in
(0:_{G(\mathbf m)}\mathcal M^r)$,  by Remark \ref{1a}. Conversely,
if $a_i=b_i$, $\omega_i-(b_i+1)g_1=\omega_i-a_ig_1-g_1=\omega'_i-g_1 \notin S'$;
again by  Remark \ref{1a}, we get $\overline{t^{\omega_i}}\notin (0:_{G(\mathbf
m)}\mathcal M^r)$.
\end{rem}

\section{Buchsbaumness in the general case}

As in the Preliminaries, $(R,\mathbf m)$ is the numerical
semigroup ring associated to a $n$-generated numerical semigroup
$S$, $G(\mathbf m)$ is its associated graded ring, $\mathcal M$ is
the homogeneous maximal ideal of $G(\mathbf m)$, $M$ is the
maximal ideal of $S$ and $r$ is the reduction number of $\mathbf
m$.

For each $i$ such that $a_i>b_i$, let $l_i=\max\left\{l\ |\
\overline{t^{\omega_i+lg_1}}\in (0:_{G(\mathbf m)}\mathcal
M^r)\right\}$.

\begin{rem} We note that $l_i$ is well defined because from
Formula (\ref {form2}) we have that
$\overline{t^{\omega_i+lg_1}}\notin (0:_{G(\mathbf m)}\mathcal
M^r)$, whenever $t^{\omega_i+lg_1}\in \mathbf m^{r-1}$.
\end{rem}

By \cite[Proposition 3.5]{DMM}, if $a_i>b_i$, then $r\ge b_i+2$.

\begin{rem} We note that $l_i\le r-2-b_i$. Indeed, $\omega_i\in b_iM$
hence $\omega_i+l_ig_1\in (b_i+l_i)M$ and $b_i+l_i\le r-2$ by
Formula (\ref {form2}).
\end{rem}

\begin{lem}\label{2a}
Let $i$ and $l_i$ as above. Then $\overline{t^{\omega_i+lg_1}}\in
(0:_{G(\mathbf m)}\mathcal M^r)$ for every $l=0,\dots,l_i$.
\end{lem}
\noindent {\bf Proof}. By hypothesis
$\overline{t^{\omega_i+l_ig_1}}\in (0:_{G(\mathbf m)}\mathcal
M^r)$ therefore, by Remark \ref{1a}, if $\omega_i+l_ig_1\in
hM\setminus (h+1)M$ then $\omega_i+l_ig_1-(h+1)g_1\in S'$. Let
$l=l_i-1$ and let us suppose $\omega_i+lg_1\in nM\setminus
(n+1)M$. If $\overline{t^{\omega_i+lg_1}}\notin (0:_{G(\mathbf
m)}\mathcal M^r)$ then $\omega_i+lg_1-(n+1)g_1\notin S'$; it
follows that $\omega_i'>\omega_i+lg_1-(n+1)g_1$. Since $n<h$ and
$l=l_i-1$, we get $\omega_i+l_ig_1-(h+1)g_1=\omega_i+lg_1-hg_1\leq
\omega_i+lg_1-(n+1)g_1<\omega_i'$. Hence
$\omega_i+l_ig_1-(h+1)g_1\notin S'$; contradiction. Using a
decreasing induction we get the thesis. \hfill $\Box$

\begin{lem}\label{2b}
The only monomials in $(0:_{G(\mathbf m)}\mathcal M^r)$ are
of the form $\overline{t^{\omega_i+lg_1}}$, with $i$ such that $a_i>b_i$.
\end{lem}
\noindent {\bf Proof}. Let $\overline{t^{c}}\in (0:_{G(\mathbf
m)}\mathcal M^r)$. Then
$\overline{t^{c}}\in\frac{(t^{g_1})^{h+1}R'\cap \mathbf
m^{h}}{\mathbf m^{h+1}}$, with $h$ such that $t^c\in \mathbf
m^h\setminus\mathbf m^{h+1}$. In particular $c\in S$, hence
$c=\omega_i+lg_1$, for some index $i$.

Let us show that the case $a_{i}=b_{i}$ is not possible. Since
$\overline{t^{\omega_{i}+lg_1}}\in (0:_{G(\mathbf m)}\mathcal
M^r)$, by Remark \ref{1a} it follows that
$\omega_i+lg_1=\alpha+(h+1)g_1\in hM\setminus (h+1)M$ with
$\alpha=\omega_i'+\mu g_1\in S'$, $\mu\ge 0$, that is
$\omega_i+lg_1= \omega_i'+(\mu+h+1) g_1\in hM\setminus (h+1)M$. In
particular, it is in $S$ and this implies $\mu+h+1\ge a_i=b_i$.
Furthermore, $\omega_i'+(\mu+h+1)
g_1=\omega_i'+b_ig_1+(\mu+h+1-b_i) g_1\in
(b_iM+(\mu+h+1-b_i)M)\setminus (h+1)M=(\mu+h+1)M\setminus (h+1)M$
and we get $\mu+h+1<h+1$. Absurd. \hfill $\Box$

\begin{cor}\label{3}
Let $G(\mathbf m)$ be not C-M. Then
$$
(0:_{G(\mathbf m)}\mathcal
M^r)=\left\langle \overline{t^{\omega_i+lg_1}}\mid a_i>b_i,\ l=0,\dots,l_i \right\rangle_k.
$$
\end{cor}
\noindent {\bf Proof}. It follows by Lemmas \ref{2a} and \ref{2b}. \hfill $\Box$

\medskip

Furthermore, by the previous corollary and by Proposition \ref{a},
we get the following characterization.

\begin{prop}\label{4a}
$G(\mathbf m)$ is Buchsbaum if and only if
$\overline{t^{\omega_i+lg_1}}\in (0:_{G(\mathbf m)}\mathcal M)$,
for every $i$ such that $a_i>b_i$ and for every $l=0,\dots,l_i$.
\end{prop}

The following proposition gives a bound for $l_i$ when $G(\mathbf
m)$ is Buchsbaum.

\begin{prop}\label{3a}
Let $G(\mathbf m)$ be Buchsbaum and $i$ such that $a_i>b_i$. Then $l_i<a_i-b_i$.
\end{prop}
\noindent {\bf Proof}. By definitions of $a_i$ and $b_i$ we have
that $\omega_i=\omega'_i+(a_i-b_i-1)g_1+(b_i+1)g_1\in
b_iM\setminus (b_i+1)M$ with $\omega'_i+(a_i-b_i-1)g_1\in S'$. By
hypothesis $(0:_{G(\mathbf m)}\mathcal M)=(0:_{G(\mathbf
m)}\mathcal M^r)$, hence $\omega_i+l_ig_1\in hM\setminus (h+1)M$
with $h\ge b_i+2l_i$. By definition of $l_i$ and by Remark
\ref{1a}, we have that $\omega_i+l_ig_1-(h+1)g_1\in S'$, hence
$\omega_i+l_ig_1-(h+1)g_1\ge \omega'_i=\omega_i-a_ig_1$ that is
$l_ig_1-(h+1)g_1\ge -a_ig_1$. Finally, $l_ig_1+a_ig_1\ge
(h+1)g_1\ge (b_i+2l_i+1)g_1$ implies $l_i<a_i-b_i$. \hfill $\Box$

\bigskip

If $a_i-b_i=1$, for every $i$ such that $a_i>b_i$, then we can
improve Proposition \ref{4a}.

\begin{prop}\label{3aa}
If $a_i-b_i=1$ for every $i$ such that $a_i>b_i$, then
\begin{center}
$G(\mathbf m)$ is Buchsbaum if and only if
$\overline{t^{\omega_i}}\in (0:_{G(\mathbf m)}\mathcal M)$ for
every such index $i$.
\end{center}
\end{prop}
\noindent {\bf Proof}. By Proposition \ref{4a}, we need only to
prove the sufficient condition. By hypothesis and by definition of
$b_i$ and $a_i$, we have that $\omega_i+g_1\in hM\setminus (h+1)M$
with $h\ge b_i+2$ and $\omega_i+g_1-(h+1)g_1=\omega_i-hg_1\notin
S'$. Hence $\overline{t^{\omega_i+g_1}}\notin (0:_{G(\mathbf
m)}\mathcal M^r)$ by Remark \ref{1a} and, by Lemma \ref{2a}, we
get $\overline{t^{\omega_i+lg_1}}\notin (0:_{G(\mathbf m)}\mathcal
M^r)$ for every $l\ge 1$. Finally, by Proposition \ref{3a} we get
the thesis. \hfill $\Box$

\begin{rem}\label{b}
We note that the last proposition does not hold if there exists an
$i$ such that $a_i-b_i\ge 2$. Indeed, let $S=\left\langle
12,19,29,104\right\rangle$. The reduction number of $S$ is $r=8$.
The only index for which $a_i>b_i$ is $i=8$ with $a_8=4$ and
$b_8=1$; moreover $\omega_8=g_4=104$. Since $\omega_8+g_j\in 3M$
for each $g_j=12,19,29,104$, then $\overline{t^{\omega_8}}\in
(0:_{G(\mathbf m)}\mathcal M)$. Anyway $G(\mathbf m)$ is not
Buchsbaum. Indeed $\overline{t^{\omega_8+g_1}}\notin
(0:_{G(\mathbf m)}\mathcal M)$ as $\omega_8+g_1=116\in 4M\setminus
5M$ and $116+g_1=128\in 5M$. On the other hand
$\overline{t^{\omega_8+g_1}}\in (0:_{G(\mathbf m)}\mathcal M^8)$
as $116+(8M\setminus 9M)\subseteq 13M$.
\end{rem}

\begin{ex}\label{5a}
Let $R$ be the semigroup ring associated to the numerical
semigroup $S=\langle 17,18,21,28,29,32,33\rangle$. We use
Proposition \ref{3aa} in order to show that $G(\mathbf m)$ is
Buchsbaum. The only indexes $i$ such that $a_i>b_i$ are $i=7,10$
and in both cases we have that $a_i=3>2=b_i$. We need to check
that $\omega_7+g_j=58+g_j\in (b_7+2)M=4M$ and
$\omega_{10}+g_j=61+g_j \in (b_{10}+2)M=4M$,  for each
$g_j=17,18,21,28,29,32,33$. Since $4M=\left\{68,\longrightarrow
\right\}$, this is clearly true and $G(\mathbf m)$ is Buchsbaum.
\end{ex}

Our next goal is to relate the Buchsbaumness of $G(\mathbf m)$ to
the length $\lambda=\lambda(H^0_{\mathcal M})$ of the $G(\mathbf
m)$-module $H^0_{\mathcal M}=(0:_{G(\mathbf m)}\mathcal M^r)$ (see
the next Proposition \ref{10}).

\begin{lem}\label{6}
We have $\lambda=1$ if and only if $(0:_{G(\mathbf m)}\mathcal M^r)=G(\mathbf m)x$, with $x\in(0:_{G(\mathbf
m)}\mathcal M)$.
\end{lem}
\noindent {\bf Proof}. Let $\lambda=1$. Clearly $N:=(0:_{G(\mathbf
m)}\mathcal M^r)$ must be principal as a $G(\mathbf m)$-module. If
$x\notin (0:_{G(\mathbf m)}\mathcal M)$, then $(0)\subsetneq
\mathcal MN$. Moreover, by the graded version of Nakayama's Lemma,
$\mathcal MN\subsetneq N$. A contradiction to $\lambda=1$.

Conversely, let $N:=(0:_{G(\mathbf m)}\mathcal M^r)=G(\mathbf
m)x$, with $x\in(0:_{G(\mathbf m)}\mathcal M)$ and let us suppose
$(0)\subsetneq H\subsetneq N$, with $H$ submodule of $N$. Then,
every $\bar h\in H$, $\bar h\ne 0$, is of the form $\bar h =\bar
gx$ with $\bar g\in G(\mathbf m)$. Since $\bar h\ne 0$, we have
that $\bar g\notin\mathcal M$, then $\bar g=\bar x_0+\bar y\in
R/\mathbf m\oplus\mathcal M$ with $x_0\ne 0$. Finally $\bar h=\bar
gx=(\bar x_0+\bar y)x=\bar x_0x+\bar yx=\bar x_0x$ and, since
$\bar x_0$ is a unit in $G(\mathbf m)$, we get $x\in H$. By the
choice of $x$, we get $N\subseteq H$; a contradiction. \hfill $\Box$

\begin{cor}\label{7}
If $\lambda\le 1$, then $G(\mathbf m)$ is Buchsbaum.
\end{cor}
\noindent {\bf Proof}. If $\lambda=0$, then $(0:_{G(\mathbf
m)}\mathcal M^r)=(0)$, that is $G(\mathbf m)$ is C-M. If
$\lambda=1$, then, by Lemma \ref{6},  $(0:_{G(\mathbf m)}\mathcal
M^r)=G(\mathbf m)x$, with $x\in(0:_{G(\mathbf m)}\mathcal M)$.
This implies $(0:_{G(\mathbf m)}\mathcal M^r)\subseteq
(0:_{G(\mathbf m)}\mathcal M)$. \hfill $\Box$

\begin{rem}\label{8}
We note that the converse of Proposition \ref{7} does not hold in
general. Indeed, let $R$ be as in the Example \ref{5a}. We showed
that $G(\mathbf m)$ is Buchsbaum. Moreover, by Corollary \ref{3}
and Proposition \ref{3a}, we have that $(0:_{G(\mathbf m)}\mathcal
M^r)=\left\langle
\overline{t^{58}},\overline{t^{61}}\right\rangle_{G(\mathbf
m)/\mathcal M}$, hence $(0)\subsetneq G(\mathbf
m)\overline{t^{58}}\subsetneq (0:_{G(\mathbf m)}\mathcal M^r)$.
\end{rem}

It is possible to relate $\lambda$ with the $l_i$'s when
$G(\mathbf m)$ is Buchsbaum.

\begin{lem}\label{9}
Let $G(\mathbf m)$ be Buchsbaum. Let $i,j$ be such that $a_i>b_i$
and $a_j>b_j$; then $\overline{t^{\omega_i+lg_1}}\in G(\mathbf
m)\overline{t^{\omega_j}}$ if and only if $i=j$ and $l=0$.
\end{lem}

\noindent {\bf Proof}. If $\overline{t^{\omega_i+lg_1}}\in
G(\mathbf m)\overline{t^{\omega_j}}$, then
$\overline{t^{\omega_i+lg_1}}=\overline u\overline{t^{\omega_j}}$.
Since $\overline{t^{\omega_j}}\in (0:_{G(\mathbf m)}\mathcal M)$
and $\overline{t^{\omega_i+lg_1}}\ne \overline {0}$, then
$\overline u\notin\mathcal M$ and $\overline u= u+\mathbf m$ with
$u\in k$; hence $\overline{t^{\omega_i+lg_1}}=\overline
u\overline{t^{\omega_j}}=\overline{ut^{\omega_j}}$. The last
equality is equivalent to the fact that $t^{\omega_i+lg_1},
ut^{\omega_j}\in \mathbf m^{b_j}\setminus\mathbf m^{b_{j+1}}$ and
$t^{\omega_i+lg_1}-ut^{\omega_j}\in \mathbf m^{b_{j+1}}$, that is
$\omega_i+lg_1=\omega_j$, but this is true only for $l=0$ and
$i=j$. \hfill $\Box$

\bigskip

\noindent The next proposition immediately follows by Corollary \ref{3} and
Lemma \ref{9}.

\begin{prop}\label{10}
If  $G(\mathbf m)$ is Buchsbaum, then $\lambda=\sum_{i\in I}(l_i+1)$ with $I=\left\{i\ |\ a_i>b_i\right\}$.
\end{prop}

\begin{cor}\label{10a}
If  $G(\mathbf m)$ is Buchsbaum, then $\lambda\le\sum_{i\in
I}(a_i-b_i)$ with $I=\left\{i\ |\ a_i>b_i\right\}$. Moreover, if
$a_i=b_i+1$ for every $i\in I$, then $\lambda=\left|I\right|$
\end{cor}
\noindent {\bf Proof}. It follows by Propositions \ref{10},
\ref{3a} and \ref{3aa}. \hfill $\Box$

\medskip
Our next aim is to study for which elements $\omega_i$ of the
Apery set of $S$, it is possible to have $a_i>b_i$; we get a
necessary condition, in the case $G(\mathbf m)$ is Buchsbaum.

Let $s\in S$ and define $\ord (s):=h$ if $s\in hM\setminus
(h+1)M$. In particular we have $\ord(\omega_i)=b_i$. We now
introduce a partial ordering on $S$ as in \cite{B}: given $u,u'\in
S$, we say that $u\le_Mu'$ if $u+s=u'$ (hence $u\le_Su'$) and
$\ord (u)+\ord (s)=\ord(u')$ for some $s\in S$.

\begin{rem}
The set of maximal elements of $\ap_{g_1}(S)$ with this partial
orde\-ring is denoted with $\map (S)$. We note that the set of
maximal elements in $\ap_{g_1}(S)$ with the usual ordering $\le_S$
is contained in $\map(S)$ and the inclusion can be strict. For
example, let $S=\langle 8,9,15\rangle$. The only maximal ele\-ment
in $\ap_{g_1}(S)$ with respect to $\le_S$ is $45$. Anyway
$\map(S)=\left\{30,45\right\}$. Note that
$\ord(45)=5>3=\ord(30)+\ord(15)$.
\end{rem}

We say that $\omega_i$ and $\omega_j$ are {\em comparable} if
$\omega_i\le_M\omega_j$ or vice versa.

\begin{rem}\label{11}
Let $G(\mathbf m)$ be Buchsbaum. If $a_i>b_i$ and $a_j>b_j$, then
$\omega_i$ and $\omega_j$ are not comparable. Indeed, if there
exists $s\in S\setminus\left\{0\right\}$ such that
$\omega_j=\omega_i+s\in S$ and $b_i+\ord (s)=b_j$, then
$\overline{t^{\omega_i}}\notin(0:_{G(\mathbf m)}\mathcal M^{\ord
(s)})=(0:_{G(\mathbf m)}\mathcal M^r)$.  A contradiction to Remark
\ref{1}.
\end{rem}

\begin{prop}\label{14}
Let $G(\mathbf m)$ be Buchsbaum. Then $a_i>b_i$ implies $\omega_i\in\map(S)$.
\end{prop}
\noindent {\bf Proof}. If $\omega_i\notin\map(S)$, then there
exists $\omega_j$ such that $\omega_i<_M\omega_j$, that is there
exists an element $s\in S\setminus\left\{0\right\}$ such that
$\omega_j=\omega_i+s$ and $b_j=b_i+\ord(s)$. By hypothesis
$\overline{t^{\omega_i}}\in (0:_{G(\mathbf m)}\mathcal M)$, hence
$\overline{t^{\omega_i}}\cdot\overline{t^{s}}=\overline{0}$; this
implies $\omega_j=\omega_i+s\in(b_i+\ord(s)+1)M=(b_j+1)M$. Absurd. \hfill $\Box$

\begin{rem} We note that the converse of the last proposition does not
hold in general. Indeed, let $R$ be the semigroup ring associated
to $S=\langle 12,19,29,104\rangle$. In this case the unique index
$i$ such that $a_i>b_i$ is $i=8$ and $\omega_8=104\in\map(S)$; but
$G(\mathbf m)$ is not Buchsbaum as showed in Remark \ref{b}.
\end{rem}

We end this section with a general remark that will be useful for
the next sections.

\begin{rem}\label{20}
Let $a_i=b_i$ and let $\omega_i=\alpha_2g_2+\cdots+\alpha_ng_n$
with $\sum_{k=2}^n\alpha_k=b_i$. By definition of $a_i$ and by the
equality $a_i=b_i$, we have that
$\alpha_2(g_2-g_1)+\cdots+\alpha_n(g_n-g_1)=\omega_i'\in\ap_{g_1}(S')$.

On the other hand, let $a_i>b_i$ and $\omega_i=\alpha_2g_2+\cdots+\alpha_ng_n$ with
$\sum_{k=2}^n\alpha_k=b_i$. In this case, $\alpha_2(g_2-g_1)+\cdots+\alpha_n(g_n-g_1)\notin\ap_{g_1}(S')$.
Indeed, $\alpha_2g_2+\cdots+\alpha_ng_n-(\sum_{k=2}^n\alpha_k)g_1=
\alpha_2g_2+\cdots+\alpha_ng_n-b_ig_1>\alpha_2g_2+\cdots+\alpha_ng_n-a_ig_1=\omega_i'$, hence
$\alpha_2(g_2-g_1)+\cdots+\alpha_n(g_n-g_1)-g_1\in S'$.
\end{rem}

\section{The $3$-generated case}
In this section we will apply and deepen our results, when the
semigroup $S$ is $3$-generated. As a by-product, we will give a
positive answer to two conjectures raised by Sapko in \cite{S}.
These two conjectures are also proved by Shen in \cite{SH} using
completely different methods.

\bigskip

Let us fix our notation for this section: $S=\langle
g_1,g_2,g_3\rangle$, with $g_1 < g_2 < g_3$; the elements in
$\ap_{g_1}(S)$ are of the form $\omega_i=hg_2+kg_3$ (with $h,k \in
\mathbb N$).

With the symbol $x\equiv y$ we will always mean that $x$ is
congruent $y$ modulo $g_1$; moreover, $x\lhd y$ (respectively
$x\rhd y$) will always mean that $x\equiv y$ and that $x<y$ (resp.
$x>y$). Finally $x\unlhd y$ (resp. $x\unrhd y$) will mean that
$x\equiv y$ and that $x\le y$ (resp. $x\ge y$).

If an element $\omega_i$ has more than one representation as a
combination of $g_2$ and $g_3$, then the representation
$hg_2+kg_3$, where $h$ is maximum, has the property that $h+k=b_i$
(this is not true if $S$ has more than $3$ generators).

We are ready to prove the main result of this section.

\begin{thm}\label{42}
Assume that $S$ is a $3$ generated numerical semigroup and assume
that $G(\mathbf m)$ is Buchsbaum. If \ $\omega_i=hg_2+kg_3$ \ is
an element of $\ap_{g_1}(S)$ such that $a_i > b_i$, then $h=0$.

In particular, there is at most one element $\omega_i \in \ap_{g_1}(S)$ such that $a_i > b_i$.
\end{thm}
\noindent {\bf Proof}. Since $G(\mathbf m)$ is Buchsbaum, by
Proposition \ref{14} we have that, if $a_i>b_i$, then
$\omega_i=hg_2+kg_3\in\map(S)$, for any representation of
$\omega_i$ as a combination of $g_2$ and $g_3$. In particular, we
consider the representation with $h$ maximum, that is $h+k=b_i$.
If we prove that $h=0$, then we get the first part of the theorem.

Assume, by contradiction, that $h>0$, hence
$(h-1)g_2+kg_3=\omega_j \in \ap_{g_1}(S)$. We note that
$\omega_j\le_M\omega_i$ as $\omega_i=\omega_j+g_2$ and $b_i=b_j+1$
(clearly $b_i\ge b_j+1$ and, if $b_i>b_j+1$, then $h+k>b_j+1\ge
(h-1)+k+1=h+k$).

Since $\omega_j\notin\map(S)$ and $G(\mathbf m)$ is Buchsbaum, we
have $a_j = b_j$ again by Proposition \ref{14}. By Remark \ref{20}
we have $$\omega'_j=(h-1)(g_2-g_1)+k(g_3-g_1).$$ By the same
remark we also get $\omega'_i \lhd h(g_2-g_1)+k(g_3-g_1)$;
recalling that $\omega'_i \in \ap_{g_1}(S')$ and $S'=\langle
g_1,g_2-g_1, g_3-g_1 \rangle$, we obtain
$\omega'_i=x(g_2-g_1)+y(g_3-g_1)$, for some nonnegative integers
$x$ and $y$. We collect this observation in the following formula
\begin{equation}
\omega'_i=x(g_2-g_1)+y(g_3-g_1)\lhd h(g_2-g_1)+k(g_3-g_1)
\end{equation}
Moreover, by definition of $a_i$ we obtain:
$$
xg_2+yg_3-(x+y)g_1+a_ig_1=\omega_i \in  \ap_{g_1}(S) ;
$$
it follows immediately that $x+y \geq a_i$. Hence, since $a_i>
b_i=h+k$, we get $x+y>h+k$.

Now, if $x\leq h$, then $y>k$ and it would follow that $\omega'_i=x(g_2-g_1)+y(g_3-g_1)>h(g_2-g_1)+k(g_3-g_1)$, in
contradiction to (4.1). Thus $x>h(>0)$ and, in particular, $x>0$.

It follows that $(x-1)(g_2-g_1)+y(g_3-g_1) \in S'$. But, again by
(4.1), $(x-1)(g_2-g_1)+y(g_3-g_1)\lhd
(h-1)(g_2-g_1)+k(g_3-g_1)=\omega'_j$. But this is a contradiction,
since $\omega'_j \in \ap_{g_1}(S')$.

Hence $h=0$ and we have proved the first part of the theorem.

Let us prove the last assertion. By the first part, we have that
the only elements $\omega_i$ for which it is possible to have
$a_i>b_i$ are of the form $jg_3$ with $\ord(jg_3)=j$ and the set
of this kind of elements is $\left\{0,g_3,\dots,kg_3\right\}$, for
some $k>0$. Since this set is a subchain of $(\ap(S),\le_M)$,
there is at most one maximal element. The thesis follows by
Proposition \ref{14}. \hfill $\Box$

\begin{rem}
The integer $k$ defined in the last part of the previous proof can
be also defined in terms of the Apery set of $S$ in the following
way:
$$
k=min\{j|\ g_2\ \mbox{divides}\ (j+1)g_3\ \mbox{or}\ (j+1)g_3-g_1 \in S\}.
$$
\end{rem}

\medskip

\noindent As a consequence of the previous theorem we obtain a
positive answer to two conjectures stated in \cite{S}, that we
collect in the following statement.

\begin{prop}\label{43}
Let $S$ be a $3$-generated numerical semigroup. Then the following condition are equivalent:
\begin{description}
  \item[(i)] $G(\mathbf m)$ is Buschsabum not C-M;
  \item[(ii)] $(0:_{G(\mathbf m)}\mathcal M^r)=G(\mathbf
  m)\overline{t^{kg_3}}$
      for some $k\ge 1$, with $\overline{t^{kg_3}}\in (0:_{G(\mathbf m)}\mathcal M)$;
  \item[(iii)] $\lambda = 1$.
\end{description}
\end{prop}

{\bf Proof}. The implication (ii) $\Rightarrow$ (iii) is obvious
by Lemma \ref{6} and the implication (iii) $\Rightarrow$ (i) is
straightforward by Proposition \ref{7}. Let us prove the
implication (i) $\Rightarrow$ (ii).

By Theorem \ref{42}, we know that there exists a unique $\omega_i$
such that $a_i>b_i$ and it is $\omega_i=kg_3$, with $k=min\{j|\
g_2\ \mbox{divides}\ (j+1)g_3\ \mbox{or}\ (j+1)g_3-g_1 \in S\}$.
Hence $\omega_i=kg_3$ is the only element in the Apery set of $S$
such that $\overline{t^{\omega_i}}\in (0:_{G(\mathbf m)}\mathcal
M)=(0:_{G(\mathbf m)}\mathcal M^r)$.

By Corollary \ref{4a} and Lemma \ref{9}, we need to prove that
$\overline{t^{\omega_i+lg_1}}\notin (0:_{G(\mathbf m)}\mathcal
M)$, for every $l\ge 1$. By Lemma \ref{2a}, it is enough to prove
it for $l=1$.

We note that, by definition of $k$, if $(k+1)g_3\in\ap(S)$, then
$g_2$ divides $(k+1)g_3$; hence $kg_3+g_3=qg_2$, with $q>k+1$.
Moreover, since $\overline{t^{kg_3}}\in (0:_{G(\mathbf m)}\mathcal
M)$ and since $kg_3\in\ap(S)$, then $kg_3+g_1=\alpha g_2+\beta
g_3$, with $\ord(kg_3+g_1)=\alpha+\beta>k+1$ (and $\beta<k$). This
implies that $g_3-g_1=(q-\alpha)g_2-\beta g_3$, i.e.
$(\beta+1)g_3=(q-\alpha)g_2+g_1$, with $\beta+1\le k$.
Contradiction against the assumption $(k+1)g_3\in\ap(S)$.

Hence we can assume $(k+1)g_3\notin\ap(S)$ and so there exists an
integer $q\ge 0$, such that $kg_3+qg_2$ is maximal in the Apery
set of $S$. Now, if $kg_3+qg_2=ug_3+vg_2$ with $v>q$, then
$(k-u)g_3=(v-q)g_2$ and this is a contradiction against the
definition of $k$. Hence $\ord(kg_3+qg_2)=k+q$ and necessarily
$q=0$ (if not $kg_3\notin\map(S)$).

In order to show that $\overline{t^{\omega_i+g_1}}\notin
(0:_{G(\mathbf m)}\mathcal M)$, it is enough to prove that
$kg_3+2g_1\notin (\alpha+\beta+2)M$ where, as above,
$kg_3+g_1=\alpha g_2+\beta g_3$ and
$\ord(kg_3+g_1)=\alpha+\beta>k+1$ (with $\beta<k$, $\alpha>0$).

If $kg_3+2g_1=ag_1+bg_2+cg_3$ with $\ord(kg_3+2g_1)=a+b+c\ge
\alpha+\beta+2$, then, by definition of $k$ and by
$kg_3\in\ap(S)$, we have $a<2$. The case $a=1$ is not possible, as
we would have $kg_3+g_1=bg_2+cg_3$ with
$b+c>\alpha+\beta=\ord(kg_3+g_1)$. Absurd. Hence $a=0$. If $c\ge
k$, then $2g_1\ge bg_2$ and so $b=1$; but this is not possible
since the case $c=k$ would give us $2g_1=g_2$, and the case $c>k$
would give us $2g_1=g_2+(c-k)g_3$. Hence $c<k$ (and $b>1$).

Since $a_i>b_i=k$, we have that $kg_3-kg_1\notin\ap(S')$, hence
$k(g_3-g_1)\rhd x(g_2-g_1)+y(g_3-g_1)\in\ap(S')$, for some
integers $x$ and $y$. If $y>0$, then $(k-y)(g_3-g_1)\rhd
x(g_2-g_1)\in S'$ and this is not possible since
$(k-y)(g_3-g_1)\in\ap(S')$, by Remark \ref{20}. Hence $y=0$ and
$k(g_3-g_1)\rhd x(g_2-g_1)\in\ap(S')$ and so there exists a $z>0$
such that $k(g_3-g_1)=x(g_2-g_1)+zg_1$, that is $\alpha g_2+\beta
g_3=kg_3+g_1=xg_2-(x-k-z-1)g_1$. Hence $xg_2=kg_3+\mu g_1=\alpha
g_2+\beta g_3+(\mu-1)g_1$ with $\mu>0$ and $\beta
g_3+(\mu-1)g_1=(x-\alpha)g_2$. Now $(x-\alpha)g_2\in\ap(S)$, since
$(x-1)g_2\in\ap(S)$; the last assertion follows by Remark \ref{20}
and by Theorem \ref{42}: the map
$$
\varphi: \ap(S)\setminus\left\{kg_3\right\}\longrightarrow \ap(S')\setminus\left\{x(g_2-g_1)\right\}
$$
defined by $\varphi (\gamma g_2+\delta g_3)=\gamma
(g_2-g_1)+\delta (g_3-g_1)$ is bijective. Since
$x(g_2-g_1)\in \ap(S')$, also $(x-1)(g_2-g_1)\in \ap(S')$; the bijection
implies that $(x-1)g_2\in\ap(S)$.

By $\beta g_3+(\mu-1)g_1=(x-\alpha)g_2\in\ap(S)$, we have that
$\mu=1$; moreover, since $\beta < k$, we get $x=\alpha$ and
$\beta=0$.

Hence $kg_3+g_1=\alpha g_2$ with $\alpha>k+1$, and
$kg_3+2g_1=bg_2+cg_3$ with $b+c>\alpha+2$, $c<k$ and $b>1$; hence
$g_1=(b-\alpha)g_2+cg_3$, so, necessarily, $c\ne 0$ and
$b<\alpha$; but this implies $g_1+(\alpha-b)g_2=cg_3$ and this is
absurd by definition of $k$ and by $c<k$. \hfill $\Box$

\bigskip

By the proof of the previous proposition, it is straightforward
that the integer $k$ of the statement (point (ii)) is the same
integer defined in Remark 4.2, hence it is determined in terms of
the Apery set of $S$:

\begin{cor}\label{44}
Let $S$ be a $3$-generated numerical semigroup. If $G(\mathbf m)$
is Buchsbaum not C-M, then $(0:_{G(\mathbf m)}\mathcal
M^r)=G(\mathbf m)\overline{t^{kg_3}}$, where the integer $k$ is
determined as follows:
$$
k=min\{j|\ g_2\ \mbox{divides}\ (j+1)g_3\ \mbox{or}\ (j+1)g_3-g_1 \in S\}.
$$
\end{cor}

Using Theorem \ref{42} we can also prove, in the case of $3$
generators, that, if $R$ is Gorenstein, then $G(\mathbf m)$ is C-M
if and only if it is Buchsbaum. This fact is also proved in
\cite{SH} using different methods.

\begin{cor}\label{45}
Let $S$ be a $3$-generated symmetric numerical semigroup. If
$G(\mathbf m)$ is Buchsbaum, then it is C-M.
\end{cor}

\noindent {\bf Proof}. By the proof of Theorem \ref{42} (and since
$G(\mathbf m)$ is Buchsbaum), it is possible to have $a_i>b_i$ only for
$\omega_i=kg_3$, with $k=min\{j|\ g_2 \ \mbox{does not divide}\\
(j+1)g_3 \ \mbox{or} \ (j+1)g_3-g_1 \in S\}$ (we underline that
$b_i=\ord(kg_3)=k$). So, by Theorem \ref{0}, we only need to show that
$a_i=b_i$.

Since $S$ is symmetric, there exists a unique maximal element
$g+g_1$ in the Apery set of $S$ (with the partial ordering $\le_S$
as in the Preliminaries). Assume that $kg_3+g_3\notin \ap(S)$; it
follows that $g+g_1=qg_2+kg_3$. Moreover, this representation is
unique as, if $qg_2+kg_3=ug_2+vg_3$, then $u>q$ and $v<k$ (if not
$v>k$ and $ug_2+vg_3\notin\ap(S)$), and we get $(k-v)g_3=(u-q)g_2$
that implies $\ord(kg_3)>k$. The uniqueness of the representation
implies that $\ord(qg_2+kg_3)=q+k$. It follows that
$kg_3\le_Mqg_2+kg_3$ and so $kg_3\notin\map(S)$, unless $q=0$. But
$g_2\in\ap(S)$ and $g_2\le_S g+g_1$, hence, if $q=0$, then
$\ord(kg_3)>k$. Hence $q \neq 0$, $kg_3\notin\map(S)$ and, by
Proposition \ref{14}, $a_i=b_i$.

Assume, now, that $kg_3+g_3\in \ap(S)$. By definition of $k$, it
follows that $(*) \ kg_3+g_3=ug_2$. Let us suppose that
$\ord(\omega_i+g_1)>b_i+1=k+1$, that is $kg_3+g_1=ag_1+bg_2+cg_3$,
with $a+b+c>k+1$. Since $\omega_i\in\ap(S)$, we get $a\le 1$ and
$c<k$. The case $a=1$ is not possible as $\ord(kg_3)=k$, hence
$a=0$, that is $(**) \ kg_3+g_1=bg_2+cg_3$ (with $b+c>k+1$). By
$(*)$ and $(**)$ it follows that $g_3-g_1=(u-b)g_2-cg_3$, i.e.
$(c+1)g_3=(u-b)g_2+g_1$. Since $c<k$, this is a contradiction
against the definition of $k$, therefore
$\ord(\omega_i+g_1)=b_i+1=k+1$. Thus
$\overline{t^{\omega_i}}\notin (0:_{G(\mathbf m)}\mathcal
M)=(0:_{G(\mathbf m)}\mathcal M^r)$ and, by Remark \ref{1},
$a_i=b_i$. \hfill $\Box$

\section{Cohen-Macaulayness and Gorensteinness in the general case}

In \cite[Theorem 2.6]{BF} the authors proved, in particular, that
$G(\mathbf m)$ is C-M iff $a_i=b_i$ for every $i=0,\dots,g_1-1$.
As a consequence of this, there is an algorithm to check whether
$G(\mathbf m)$ is C-M or not:

1) compute $hM$ for $h=1,\dots,r$,

2) find $\ap_{g_1}(S)$ and $\ap_{g_1}(S')$,

3) determine $a_i$ and $b_i$ for $i=0,\dots,g_1-1$

4) compare $a_i$ and $b_i$. If there exists $i$ such that
$a_i>b_i$ then $G(\mathbf m)$ is not C-M. If not, it is C-M.

\medskip

\noindent In the next proposition we improve the characterization
and the algorithm above, showing that in 3) it is sufficient to
determine $a_i$ and $b_i$, just for those $i$ such that
$\omega_i\in\map(S)$.

\begin{prop}\label{21}
$G(\mathbf m)$ is C-M if and only if $a_i=b_i$, for those $i$ such that $\omega_i\in\map(S)$.
\end{prop}
\noindent {\bf Proof}. By Theorem \ref{0}, we have only to
prove the sufficient condition. Assume that $a_i=b_i$, for those
$i$ such that $\omega_i\in\map(S)$ and let
$\omega_j=\alpha_2g_2+\cdots+\alpha_ng_n\notin\map(S)$. Then there
exists $\omega_i=\beta_2g_2+\cdots+\beta_ng_n\in \map(S)$, with
$\sum_{k=2}^n\beta_k=b_i=a_i$, such that
$\omega_j+\eta_2g_2+\cdots+\eta_ng_n=\omega_i$ and
$b_j+\ord(\eta_2g_2+\cdots+\eta_ng_n)=b_i$.

By Remark \ref{20}, if $a_j>b_j$, \
$\alpha_2(g_2-g_1)+\cdots+\alpha_n(g_n-g_1)\notin\ap_{g_1}(S')$; \ on the other hand, \
$a_i=b_i$ implies that \ $\beta_2(g_2-g_1)+\cdots+\beta_n(g_n-g_1)=\omega_i'\in\ap_{g_1}(S')$.

Finally, by
$\sum_{k=2}^n(\alpha_k+\eta_k)g_k=\sum_{k=2}^n\beta_kg_k$ and
$\sum_{k=2}^n\alpha_k+\sum_{k=2}^n\eta_k=\sum_{k=2}^n\beta_k$, we
get
$\sum_{k=2}^n\alpha_k(g_k-g_1)+\sum_{k=2}^n\eta_k(g_k-g_1)=\sum_{k=2}^n\beta_k(g_k-g_1)=\omega_i'\in\ap_{g_1}(S')$.
Contradiction.  \hfill $\Box$

\begin{ex}
Let $S=\langle 10,13,14\rangle$. Then $G(\mathbf m)$ is C-M as
$\map(S)=\left\{\omega_5=55,\omega_9=39\right\}$ and $a_5=b_5=4$
and $a_9=b_9=3$.
\end{ex}

\begin{rem}
Proposition \ref{21} does not hold in general if we only consider
the maximal elements in the Apery set of $S$ w.r.t. $\le_S$. Let
$R$ be the semigroup ring associated to  $S=\langle
7,8,9,19\rangle$. The maximal elements in the Apery set of $S$ are
$\left\{\omega_3=17, \omega_6=27\right\}$ and $a_3=b_3=2$ and
$a_6=b_6=3$. Anyway  $G(\mathbf m)$ is not C-M as $\omega_5=19$
and $a_5=2>1=b_5$.
\end{rem}

In Corollary \ref{45}, we showed that in the $3$-generated case,
if $R$ is Gorenstein, then the properties for $G(\mathbf m)$ to be
C-M and Buchsbaum are equivalent.

\begin{rem}\label{45a} We note that in the $n$-generated case,
Corollary \ref{45} is not true. Let us consider the symmetric
numerical semigroup $S=\langle 8,9,12,13,19\rangle$. The only
index $i$ for which $a_i>b_i$ is $i=3$ (in particular $G(\mathbf
m)$ is not C-M); more precisely we have $a_3=2$, $b_3=1$ and
$\omega_3=19$. Since $\overline{t^{19}}\in (0:_{G(\mathbf
m)}\mathcal M)$, then $G(\mathbf m)$ is Buchsbaum by Proposition
\ref{3aa}.
\end{rem}

Anyway, if we force the elements of $\map(S)$ to have all the same
order, then Corollary \ref{45} is true in the $n$-generated case.
A numerical semigroup $S$ is called $M$-{\em pure} if every
element in $\map(S)$ has the same order (cf. \cite{B}).
In this case it is clear that $\map(S)$ coincides with the set of the maximal elements of $\ap(S)$
with respect to $\leq_S$.

\begin{prop}\label{46}
Let $S$ be a $M$-pure symmetric numerical semigroup. If $G(\mathbf
m)$ is Buchsbaum, then it is C-M.
\end{prop}

\noindent {\bf Proof}. Let
$\ap(S)=\left\{\omega_0, \dots, \omega_{g_1-1}\right\}=\left\{0<v_1<\cdots<v_{g_1-1}=g+g_1\right\}$ (where $<$ is the
natural ordering in $\mathbb N$). Since $S$ is $M$-pure and
symmetric, we get $\map(S)=\left\{v_{g_1-1}\right\}$.
Since $G(\mathbf m)$ is Buchsbaum, it is possible to have
$a_i>b_i$ only for $\omega_i=v_{g_1-1}$. So, by Theorem \ref{0}, we only need to show that
$a_i=b_i$.

Let us consider $\omega'_i$. If it is a minimal generator of $S'$ (different from
$g_1$, because $\omega'_i\in\ap(S')$), then
$\omega'_i=g_t-g_1$, for some generator $g_t$ of $S$ ($t \neq
1$); hence $\omega_i=g_t$ and this implies
$a_i=1=b_i$.

On the other hand, if $\omega'_i$ is not a minimal generator of $S'$,
it can be written as a sum of two elements of $S'$, that are necessarily elements
of $\ap(S')$. Hence we have an equality $\omega'_i=\omega'_j+\omega'_h$,
for some $j,h \neq i$. It follows that $\omega_i=g+g_1\equiv \omega_j+\omega_h$ and
by the simmetry of $S$ we immmediately get $\omega_i=g+g_1= \omega_j+\omega_h$.
It follows that
$$
\aligned \omega'_i&+a_ig_1=\omega_i=\omega_j+\omega_h=\omega'_j+\omega'_h+(a_j+a_h)g_1= \\
&= \omega'_j+\omega'_h+(b_j+b_h)g_1=\omega'_i+(b_j+b_h)g_1 \ . \endaligned
$$
Hence $a_i=b_j+b_h \leq b_i \leq a_i$ and so $a_i=b_i$. \hfill $\Box$

\medskip

As an immediate corollary of the last proposition we can improve
\cite[Corollary 3.20]{B} in which the author proved that:
\begin{center}
$G(\mathbf m)$ is Gorenstein $\Longleftrightarrow$ $S$ is
symmetric, $M$-pure and $G(\mathbf m)$ is C-M.
\end{center}

\begin{cor}\label{47}
$G(\mathbf m)$ is Gorenstein if and only if $S$ is symmetric, $M$-pure and $G(\mathbf m)$ is Buchsbaum.
\end{cor}

Let $J$ be a parameter ideal of a Noetherian local ring $(R,
\mathbf m)$; the index of nilpotency of $\mathbf m$ with respect
to $J$ is defined to be the integer $s_J(\mathbf m)= \text{min}\{
n |\ \mathbf{m}^{n+1} \subseteq J \}$.

If $J=(t^{g_1})$, then $J$
is a reduction of $\mathbf{m}$, $s_J(\mathbf
m)=\text{max}\{ord(\omega_i) | \ \omega_i \in \ap(S)\}$ and
$s_J(\mathbf m)\leq r$, where $r$ is the reduction number of $R$
(see, e.g., \cite{B}). In \cite[Theorem 3,14]{B}, it is also
proved that, with this choice of $J$,
\begin{center}
if $S$ is $M$-pure, then $G(\mathbf m)$ is C-M \
$\Longleftrightarrow$ \ $s_J(\mathbf m)=r$.
\end{center}
Hence, combining the previous results we immediately get
\begin{center}
$S$ is $M$-pure, symmetric and $s_J(\mathbf m)=r$ \
$\Longleftrightarrow$
\\ $S$ is $M$-pure, symmetric and $G(\mathbf m)$ is Buchsbaum.
\end{center}
It is natural to ask if, in the previous equivalence, one can skip
one or both the condition $S$ symmetric and $S$ $M$-pure.

It is easy to see that the condition $G(\mathbf m)$ Buchsbaum does
not imply $s_J(\mathbf m)=r$: if $S=\langle 4,5,11 \rangle$ then
$G(\mathbf m)$ is Buchsbaum (since $r=3$; cf. \cite[Proposition
7.7]{Go2}), but $s_J(\mathbf m)=2 < r$. Also the implication
$s_J(\mathbf m)=r \Rightarrow G(\mathbf m)$ Buchsbaum is false: if
$S=\langle 9,10,11,23 \rangle$, $s_J(\mathbf m)= r=4$, but
$G(\mathbf m)$ is not Buchsbaum, as follows by Proposition
\ref{14}, since $2=a_5 > b_5=1$ and $\omega_5=23 \notin \map(S)$.

Since $S=\langle 9,10,11,23 \rangle$ is a symmetric numerical
semigroup, the same example shows that $$S \text{ symmetric and }
s_J(\mathbf m)= r \nRightarrow S \text{ symmetric  and } G(\mathbf
m) \text{ Buchsbaum}.
$$
As for the converse, we do not have counterexamples nor an
evidence that it should be false.

\medskip
\begin{qst} Let $(R,\mathbf m)$
be a numerical semigroup ring with associated semigroup $S$.
Assume that $S$ is symmetric and $G(\mathbf m)$ is Buchsbaum; is
it true that $s_J(\mathbf m) = r$?
\end{qst}

\medskip

Finally, by  \cite[Theorem 3,14]{B} we know that, if $S$ is
$M$-pure, $s_J(\mathbf m) = r$ is equivalent to $G(\mathbf m)$ C-M; hence
if $S$ is  $M$-pure and $s_J(\mathbf m) = r$, then $G(\mathbf m)$ is
Buchsbaum; conversely, if $S$ is $M$-pure and $G(\mathbf m)$ is
Buchsbaum, we get that $s_J(\mathbf m) = r$ if and only if $G(\mathbf m)$
is C-M. We do not have any example of an $M$-pure numerical
semigroup such that $G(\mathbf m)$ is Buchsbaum not C-M.

\medskip
\begin{qst} Let $(R,\mathbf m)$ be a numerical semigroup
ring with associated semigroup $S$. Assume that $S$ is $M$-pure
and $G(\mathbf m)$ is Buchsbaum; is it true that $G(\mathbf m)$ is
C-M?
\end{qst}

After the paper was accepted for publication, we recieved by Y. H.
Shen, the following affirmative answer to Question 5.8, that we
publish with his permission.

\begin{prop}\label{answer}(Shen)
Let $(R,\mathbf m)$ be a numerical semigroup
ring with $M$-pure associated semigroup $S$.
If $G(\mathbf m)$ is
Buchsbaum, then $s_J(\mathbf m) = r$.
\end{prop}

\noindent {\bf Proof.} Since S is $M$-pure, for every $\omega_i \in
\map(S)$, we have $ord(\omega_i) = s_J(\mathbf m)\leq r$. For
simplicity of notation let us denote $s_J(\mathbf m)=s$. If $s<
r$, then $\mathbf m^{s+1}\neq t^{g_1}\mathbf m^s$. Indeed, we will
have $(t^{g_1}) \supseteq \mathbf m^{s+1} \supsetneq
t^{g_1}\mathbf m^s$. Hence, there exists a monomial $x \in \mathbf
m^{s+1}$, so that $x = (t^{g_1})y$, but $y \notin \mathbf m^s$.
Thus, $\overline t^{g_1} \overline y = 0 \in G(\mathbf m)$ and,
since $(t^{g_1}) \supseteq \mathbf m^{r}$ (by $r>s$), $y \in
H^0_{\mathcal M}=(0:_{G(\mathbf m)}\mathcal M^r)$.

Now, by Lemma 3.4, $y = t^{\omega_i+lg_1}$ for some $i$, such that
$a_i > b_i$, and some $l$, with $0 \leq l \leq l_i$. Meanwhile,
$G(\mathbf m)$ is Buchsbaum, hence, by Proposition 3.19, for this
index $i$ we have $\omega_i \in \map(S)$; moreover, by $M$-purity we
get $ord(\omega_i) = s$. But then $ord(\omega_i)+l \leq ord(y) <
s$, a contradiction. Thus we have $s = r$, as expected. \hfill
$\Box$

\bigskip

We can collect the previous results and discussion in the
following corollary.

\begin{cor}\label{answer2}
Let $(R,\mathbf m)$ be a numerical semigroup
ring with $M$-pure associated semigroup $S$. Then the following
conditions are equivalent:
\begin{description}
  \item[(i)] $G(\mathbf m)$ is Buschsabum;
  \item[(ii)] $G(\mathbf m)$ is CM;
  \item[(iii)] $s_J(\mathbf m) = r$.
\end{description}
\end{cor}

\end{document}